\documentclass{article}

\newcommand{\R}{I\!\! R}

\newcommand{\Z}{I\!\! Z}

\newcommand{\Q}{I\!\! Q}

\begin{document}

Tsemo Aristide

 Visitor, University of Toronto

 100  St. George Street

Toronto, Ontario Canada

M5S 3G3

tsemoaristide@hotmail.com

\bigskip

\bigskip

\centerline{\bf Affine manifolds, lagrangian manifolds.}

\bigskip

\centerline{\bf Abstract.}

{\it Let $(M,\omega)$ be a symplectic manifold endowed with a
lagrangian foliation ${\cal L}$, it has been shown by Weinstein [16]
that the symplectic structure of $M$ defines on each leaf of ${\cal L}$,
a connection which curvature and torsion forms vanish identically.
Suppose that $L_0$ is a compact leaf which Weinstein connection is
geodesically complete, Molino and Curras-Bosch [2] have classified
germs of such lagrangian foliation around $L_0$. In this paper we
extend this classification without supposing the completness of
the compact leaf. The Weinstein connection is dual to the Bott
connection, this enables to  relate the conjecture of Auslander
and Markus to transversally properties of these foliations.}

\bigskip

\bigskip

{\bf 1. The Weinstein connection.}

\bigskip

In this section, we shall describe following Dazord [3], the Weinstein
connection defined on the leaves of a lagrangian foliation.
Let $(M,\omega)$ be a symplectic manifold endowed with a
lagrangian foliation ${\cal L}$,  $L_0$ a leaf of ${\cal L}$,
$\chi({\cal L})$ the vector fields tangent to ${\cal L}$ and
$\chi(n({\cal L}))$ the sections of the normal bundle $n({\cal L})$ of
${\cal L}$.  One can define the Bott connection $\hat\nabla:\chi({\cal
L})\otimes \chi(n({\cal L}))$
$$
\hat\nabla_XY=u[X,Y'],
$$
where $X$ is an element of $\chi({\cal L})$,  $Y'$ is an element
of $\chi(M)$ over $Y$ and $u$ the canonical projection
$\chi(M)\rightarrow \chi(n({\cal L}))$. The connection $\hat\nabla$
induces a connection $\nabla'$ on $\chi(n({\cal L})^*)$ defined by:

$$
(\nabla'_Xf)(t)=L_X(f(t))-f(\hat\nabla_Xt).
$$
    The symplectic duality allows to identify
$\chi({\cal L})$ to $\chi(n({\cal L})^*)$, since ${\cal L}$  is
lagrangian, we deduce that the connection $\nabla'$ endows each
leaf $L_0$ of  ${\cal L}$ with a connection ${\nabla}_{L_0}$. It
has be shown by Weinstein that the curvature and torsion forms of
these connections vanish. This is equivalent to saying that the
differentiable structure of each leaf of ${\cal L}$ is defined by
an atlas which coordinates change are affine maps.

Consider a darboux system of coordinates $(q_1,..,q_n,p_1,..,p_n)$
 adapted to the foliation i.e such that the foliation is defined
 by the equations $dq_1=..=dq_n=0$. The coordinates $(p_1,..,p_n)$
 defines also the requested affine structure.

 The connection $\nabla$ is the dual of the Bott connection. This
 gives an isomorphism between the holonomy of $\nabla_{L_0}$ and the
 infinitesimal holonomy of $L_0$. The purpose of this paper is to
 interpret properties and conjectures of affine manifolds theory
 in terms of the transverse geometry of lagrangian foliations  and
 to study the neighbourhood of the compact leaf whose Weinstein
 connection is complete using the reduction principle, in Molino
 [9] and ideas of geometric non abelian higher cohomology as described in
 Tsemo  [15].

 \bigskip

 {\bf 2. The Auslander conjecture and the growth.}

 \medskip

 It has been conjectured by Auslander that the fundamental group
of a compact and complete affine manifold is polycyclic. This is
equivalent to say that the linear holonomy of the complete
structure is polycyclic.

The Auslander conjecture have been studied by many authors: It has
been proven in the following case:

- The dimension of $M$ is less than $3$, by Fried,and Goldman

- The Dimension of $M$ is less than $6$ by Abels, Soifer an
Margulis

- The linear holonomy of the leaves is includes in $O(n-1,1)$ by
Goldman and Kamishima

- The linear holonomy is included in an algebraic group product of
algebraic groups of rank less than $1$ by Margulis,

- The affine automorphisms group is of codimension greater than
$1$ by Tsemo.

\medskip

\medskip

{\bf Definition 2.1.}

 Let ${\cal L}$ be a foliation on a
compact manifold $M$, and $(U_i)_{i\in I}$ a finite atlas of $M$
such that the  restriction of ${\cal L}$ to each $U_i$ is simple,
consider a leaf $L_0$ of ${\cal L}$, and  denote by $P_i$ a connected
 component the intersection of $U_i$ and $L_0$. A path in $L_0$ associated
  to $(U_i)_{i\in I}$ is a family of open set $(P_1,...,P_k)$ such that $P_i\cap P_{i+1}$ is
not empty.

 The growth map $\gamma_{U_i}$ is
the map such that $\gamma_{U_i}(n)$ is the number of $P_j$
relatively to the atlas $(U_i)_{i\in I}$ which can be joined to
$P_i$ by a path which length is less than $n$. The growth of the
$L_0$ is the growth of $\gamma_{U_i}$.

\medskip

 Suppose that $M$ is compact, and $L_0$ is a compact leaf of
 ${\cal L}$, Haefliger has shown that
  there is a local transversal $T$ to ${\cal L}$
 which intersects $L_0$, a representation
 $\rho:\pi_1(L_0)\rightarrow Diff(T)$, such that the restriction
 of ${\cal L}$ to a neighbourhood $V$ of $L_0$ is the suspension
 foliation of the quotient of $\hat L_0\times T$ by $\pi_1(L_0)$,
 where $\hat L_0$ is the universal cover of $L_0$, and the action
 of $\pi_1(L_0)$ on $\hat L_0\times T$ is given by the Deck
 transformations on $\hat L_0$, and $\rho$ on $T$.

 \medskip

 {\bf Proposition 2.1.}

 {\it Suppose that the leaf $L_0$ is compact, then the growth of a
 leaf $L_1$ in $V$ (defined in the previous paragraph) as above is
  bounded by the growth of
 $\rho(\pi_1(L_0))$.}

 \medskip

 {\bf Proof.}

The transversal can be chosen compact, we conclude using $1.29$ of
Godbillon [6].

\medskip

{\bf Proposition 2.2.}

{\it Suppose that the growth of the leaves of ${\cal L}$ in a
neighbourood of $L_0$ is polynomial, then the fundamental group of
$L_0$ is polycyclic.}

\medskip

{\bf Proof.}

Suppose that the fundamental group is not polycyclic, then the
Tits alternative applied to the linear holonomy, implies that there exists a free subgroup of
$\pi_1(L_0)$ generated by $2$ elements, this implies that the
infinitesimal holonomy and the holonomy  of the leaf $L_0$
contains a free subgroup generated by $2$ elements. the growth of
all the leaves cannot be polycyclic. see also Plante and Thurston [12].

\medskip

We can set the following conjecture which implies the Auslander conjecture:

\medskip

{\bf Conjecture 2.3.}

{\it A compact  affine manifold is complete if and only if it is a leaf of a
lagrangian foliation of a compact symplectic  manifold such that the growth of
all the leaves are polynomial.}

\medskip

{\bf Remark.}

 There exists compact affine manifolds whose fundamental groups are
 not polycyclic, examples of such affine manifolds are the product
 of surfaces of genus greater than $2$ by a circle.

\medskip

{\bf Proposition 2.4.}

{\it Let $(L_0,\nabla_{L_0})$ be a compact  affine
manifold whose linear holonomy is contained in $Gl(n,{\Z})$, then
$(L_0,\nabla_{L_0})$ is a leaf of a lagrangian foliation defined
on a compact symplectic manifold.}

\medskip

{\bf Proof.}

Let $\hat L_0$ be the universal cover of $L_0$, the cotangent
bundle $T^*L_0$ of $L_0$ is the suspension of ${\R}^n$ over
$L_0$ by the action
$$
h:\pi_1(L_0)\rightarrow Gl(n,{\R}),
$$
 $$
 \gamma\rightarrow
^t(L(h_{L_0})(\gamma^{-1})),
$$ where $L(h_{L_0})$ is the linear
holonomy of $(L_0,\nabla_{L_0})$. Since the image of $L(h_{L_0})$
is contained in $Gl(n,{\Z})$, the image of $h$ is also contained
in $Gl(n,{\Z})$. Consider now translations $t_{e_1}$,..,$t_{e_n}$
of ${\R}^n$, where $e_1$,...,$e_n$ is a basis of ${\R}^n$. The
quotient $(M,\nabla_M)$ of $T^*L_0$ fiber by fiber by the group
generated by $t_{e_1}$,...,$t_{e_n}$ is a compact affine manifold
which is a suspension of the torus over $L_0$. Let $D$ be the
developing map of $L_0$, then the canonical symplectic form of
${\R}^{2n}$ such that the affine space $x\times {\R}^n$, an the
open set of affine spaces $D(\hat L_0)\times y$ are lagrangian is
preserved by the holonomy of $(M,\nabla_M)$. This implies that the
pulls back to $\hat M$ of this symplectic form defines on
$(M,\nabla_M)$  a symplectic affine manifold, and the suspension
foliation is lagrangian.

\bigskip

Let $(M,\omega)$ be a compact symplectic manifold endowed with a
lagrangian foliation ${\cal L}$, which has a compact and complete
leaf $(L_0,\nabla_{L_0})$. Weinstein has shown that there exists a
neighbourhood $U$ of $L_0$ which is symplectomorph to a
neighbourhood of the $0$ section of the cotangent bundle of $L_0$.
The infinitesimal holonomy of the restriction of the horizontal
foliation (if $(x_1,..,x_n)$ is a system of affine coordinates of
$L_0$, then the horizontal foliation is defined by
$(dx_1=..=dx_n=0$.) of $T^*L_0$ to $U$, and the one of ${\cal L}$
coincide, we deduce that the Auslander conjecture is true if the
combinatoric growth of the horizontal foliation of $T^*L_0$ is
polynomial. To generalize proposition 2.4, one may also find cocompact
affine action of $T^*L_0$.

\medskip

Milnor [8] asked wether the fundamental group of a complete affine
manifold is polycyclic. Margulis [7] has constructed an example of a
compact and complete affine manifold which fundamental group is a
free group generated by $2$ elements. Other  examples of free
groups which are fundamental groups of complete affine manifolds
have been constructed by Goldman, Charette and Drumm [1]. One may ask
wether those manifolds are leaves of lagrangian foliations defined
on compact symplectic manifolds.

If a compact affine manifold $L_0$ is not complete, then its
fundamental group may not be polycyclic, it is the case of the
product of a surface whose genus is greater than $2$ by a circle.
Are such manifolds leaves of lagrangian foliations on compact
manifolds?

\bigskip

{\bf Examples of compact symplectic manifolds endowed with
lagrangian foliations which has a compact leaf.}

\medskip

We will define first two lagrangian foliations on the torus $T^4$
endowed with a compact leaf.

The first example is the flat riemannian structure, we consider
the quotient of ${\R}^4$ by the translations $t_{e_1},..,t_{e_4}$.
The symplectic structure is the one induced by the one  of
${\R}^4$ defined by the form $\omega_0=dx_1\wedge dx_3+dx_2\wedge
dx_4$. The lagrangian foliation is the foliation by torus which
are the projection by the cover map of the affine spaces parallel
to $vect(e_1,e_2)$.

\medskip

The second example is the quotient of ${\R}^4$ by the
transformations:

$$
h_1(x_1,x_2,x_3,x_4)=(x_1+1,x_2,x_3,x_4)
$$

$$
h_2(x_1,x_2,x_3,x_4)= (x_1+x_2,x_2+1,x_3,-x_3+x_4)
$$

$$
h_3(x_1,x_2,x_3,x_4)= (x_1,x_2,x_3,x_4+1)
$$

$$
h_4(x_1,x_2,x_3,x_4)= (x_1+x_2,x_2,x_3+1,-x_3+x_4)
$$

The symplectic structure is the projection of $\omega_0$. The
lagrangian foliation here is defined by the projection of the
affine space parallel to $vect(e_1,e_2)$. The compact leaf is the
projection of the affine space parallel to $vect(e_1,e_2)$ which
contains $0$.

\medskip

We consider the quotient of ${\R}^4$, by $h_1$, $h_2$, $t_{e_3}$
and $t_{e_4}$, it is a symplectic manifold endowed with the
projection of $\omega_0$, the quotient of the affine space
parallel to $vect(e_1,e_2)$ is the lagrangian foliation.

\bigskip

{\bf 3. Markus Conjecture and transversely measurable foliations.}

\medskip

Markus has conjectured that a compact affine manifold is complete
if and only if it is unimodular. Unimodular means that its linear
holonomy group is a subgroup of $Sl(n,{\R})$.

Suppose that the compact and unimodular affine manifold
$(L_0,\nabla_{L_0})$ is a leaf of a lagrangian foliation ${\cal
L}$ defined on the compact manifold $M$. Let $V$ be the
neighbourhood of $L_0$, on which ${\cal L}$ is the quotient of
$\hat L_0\times T$, by the action of $\pi_1(M)$, where $T$ is a
closed lagrangian transversal to ${\cal L}$, $\pi_1(M)$ acts on
$\hat L_0$ by the deck transformations, and on $T$ by the holonomy
representation of the foliation. Since the symplectic duality
defines an isomorphism between the infinitesimal holonomy of
${\cal L}$ and the linear holonomy of the affine structure of
$L_0$, we deduce that the holonomy of ${\cal L}$ in $L_0$
preserves a measure on $T$. This leads to the folowing conjecture which
implies the Markus conjecture:

\bigskip

{\bf Conjecture 3.1.}

{\it A compact affine manifold $L_0$ is complete if and only if it
is the leaf of a lagrangian foliation which admits a transverse
measures which support contains the satured set of a transversal to $L_0$.}

\medskip

It has been shown by Plante that the growth of a leaf contained in the support
of a transverse measure of a codimension $1$
foliation on a compact manifold is polynomial. This result is not true for
codimension greater than $1$ as shows certain  suspensions of $S^2$.

A complete affine manifold which fundamental group is polycyclic is unimodular,
one may ask if a compact and complete affine manifold is a leaf of a lagrangian
foliation on a compact manifold which admits a transverse measure which
support  contains a satured space of a transversal to the compact leaf then the growth of the
leaves are polynomial?

\medskip

{\bf Proposition 3.2.}

{\it Let $(M,\nabla_M,\omega_M)$ be a symplectic affine manifold,
that is an affine manifold $(M,\nabla_M)$ endowed with a
symplectic form $\omega_M$ parallel respectively to $\nabla_M$.
Suppose that there exists in $M$ an affine lagrangian foliation
which has a compact leaf $L_0$, and such that the fundamental
group of $L_0$ is normal in $\pi_1(M)$, then $M$ is the quotient
of the cotangent bundle of $L_0$ by a discrete group of affine
symplectomorphisms.}

\medskip

{\bf Proof.}

Let $x_0$ be an element of $L_0$, the Weinstein connection is the
restriction of the connection of $(M,\nabla_M)$ to $L_0$. Let
${\R}^{2n}$ be the universal covering space of $M$, endowed with
the pulls back of the symplectic form of $M$. The
Dazord theorem implies that the quotient of ${\R}^{2n}$ by
$\pi_1(L_0)$ is the cotangent bundle of $L_0$. Since we have
supposed that $\pi_1(L_0)$ is normal in $\pi_1(M)$, we conclude
that $(M,\nabla_M)$ is the quotient of $T^*L_0$ by the cocompact
group $\pi_1(M)/\pi_1(L_0)$.

\bigskip

{\bf 4. Classification of symplectic manifolds endowed with a
lagrangian foliation.}

\bigskip

Let $(M,\omega)$ be a symplectic manifold endowed with a
lagrangian foliation ${\cal L}$ which has a compact leaf $L_0$,
suppose that there exists  a parallel action  respectively to the
Weinstein connection of the torus $T^1$ on $L_0$, Molino [9] has extended
the action of $S^1$ to $M$, and
constructed a Marsden-Weinstein reduction $(M_1,\omega_1)$ endowed
with a Lagrangian foliation ${\cal L}_1$ which has a compact leaf
$L_1$ which is the basic space of an affine bundle which total
space is $L_0$.

Suppose that the leaf $L_0$ is diffeomorphic to the
$n-$dimensional torus, then  there exists a left symmetric
algebra (associative since commutative) ${\cal H}$, which defines
a left complete invariant affine structure on the commutative
group $H$ such that $L_0$ is the quotient of $H$ by a lattice
$\pi_1(L_0)$. The Lie group $H$ is also the connected component of
the group of affine automorphisms of $(L_0,\nabla_{L_0})$, under
this identification, the associative product of ${\cal H}$ if the
restriction of the associative product of $aff({\R}^n)$ which is
defined by
$$
(C,c)(D,d)=(CD,C(d))
$$
where $C,D$ are elements of $gl(n,{\R})$ and $c,d$ are elements of
${\R}^n$.

Given generators, $\gamma_1,...,\gamma_n$ of $\pi_1(L_0)$, there
exists elements $(C_1,c_1),..,(C_n,c_n)$ such that
$\gamma_i=exp((C_i,c_i))$, we will say that the associative
structure is rational if and only if the ${\Q}-$vector space
generated by $(C_1,c_1),..,(C_n,c_n)$ is stable by the associative
structure.

\medskip

{\bf Proposition 4.1.}

{\it Let $M$ be a complete affine manifold diffeomorphic to the
torus, suppose that  the associative algebra which defines its
affine structure is rational, then there exists a parallel action
of $T^1$ on $M$.}

\medskip

{\bf Proof.}

Let $(C,c)$ be an element of the ${\Q}-$vector space generated by
$(C_1,c_1),..,(C_n,c_n)$ such that $(C,c).(C,c)=(C^2,C(c))=0$. Let
$\gamma_1,...,\gamma_p$ be generators of $\pi_1(M)$ which defines the rational
structure, with $\gamma_i=exp((C_i,c_i))$. We have
$(C,c)=\alpha_1(C_1,c_1)+..+\alpha_p(C_p,c_p)$, since ${\cal H}$
is commutative, where $\alpha_i$ is a rational. There
exists an integer $p$, such that $p\alpha_i$ is an integer, this
implies that  $exp(pu)\in\pi_1(M)$, we conclude
using Tsemo [14].

\medskip

{\bf Corollary 4.2.}

{\it Let $(M,\omega)$ be a symplectic manifold endowed with a
lagrangian foliation ${\cal L}$ which has a compact leaf $L_0$
which affine structure is isomorphic to the rational affine
structure of a complete torus, then there exists a sequence
$(M_l,\omega_l)$,...,$(M_1,\omega_1)$ of symplectic manifolds such
that $(M_l,\omega_l)$ is $(M,\omega)$, $(M_1,\omega_1)$ is the two
dimensional torus, and $(M_i,\omega_i)$ is the
Marsden-Weinstein reduction of $(M_{i+1},\omega_{i+1})$ acted on
by the circle.}

\medskip

{\bf Proof.}

Let ${\cal H}$ be  the rational associative commutative
subalgebra ${\cal H}$ of $aff({\R}^n)$   which defines the
affine structure of ${\cal L}_0$.

 We have seen that there exists a translation $t_u\in \pi_1(L_0)$,
 the quotient of $L_0$ by $t_u$ is a $n-1-$dimensional torus
 $L_1$, which affine structure is defined by the quotient ${\cal H}_1$ of the
 associative algebra ${\cal H}$ by the Lie algebra of $t_u$. Since
 ${\cal H}$ is rational, ${\cal H}_1$ is rational. Then if one
 performs a Marsden-Weinstein reduction on a $M$ acted on by the
 extended action of $t_u$, he obtains a symplectic manifold $M_1$
 which has a lagrangian foliation which has a compact leaf $L_1$
 acted on by a parallel action of the circle. We can continue this
 procedure until a torus.

\bigskip

There are complete affine structures on the $2-$dimensional torus
$T^2$, which holonomy group does not contain a translation.

Consider the $2-$dimensional commutative Lie subgroup $H$ [4], of
$Aff({\R}^2)$ which elements are $f_{s,t}(x,y)=(x+sy+ {s^2\over
2}+t,y+s)$. This group acts simply transitively on ${\R}^2$. Let
$h$ be a non rational real number, the subgroup $I$ of $H$
generated by $f_{h,0}$ and $f_{1,1}$ is a lattice of $H$. The
quotient of $H$ by $I$ is a complete affine torus which holonomy
does not have a translation.

\bigskip

More generally,
Let $(L_n,\nabla_{L_n})\rightarrow
...\rightarrow(L_1,\nabla_{L_1})$ be a sequence of compact and
complete affine manifolds where
$(L_{i+1},\nabla_{L_{i+1}})\rightarrow (L_i,\nabla_{L_i})$  is an
affine bundle which fiber is the circle such that
$(L_{i+1},\nabla_{L_{i+1}})$ is a leaf of a lagrangian foliation
defined on the symplectic manifold $M_{i+1}$, and the action of
the circle on $L_{i+1}$ extends to an hamiltonian action on
$M_{i+1}$ such that the Marsden-Weinstein reduction of $M_{i+1}$
is the symplectic manifold $M_i$ endowed with a lagrangian
foliation which has as compact leaf $L_i$. In this part we will
classify the sequence germs of neighbourhood of $L_i$ in $M_i$
for  such sequence of symplectic manifolds
$M_n\rightarrow...\rightarrow M_1$ for a given sequence of affine
manifolds $(L_n,\nabla_{L_n})\rightarrow...\rightarrow
(L_1,\nabla_{L_1})$.

\bigskip

{\bf The classification of the neighbourhood of a compact leaf of
a lagrangian foliation.}

\bigskip

In this part we recall the classication did by Molino and
Curras-Bosch that we will generalize.

\medskip

Consider a compact and complete $n-$dimensional affine manifold
$L_0$, $x_0$ an element of $L_0$, and $h_{\nabla_{L_0}}$ the
linear holonomy of $L_0$,  we
identify $T_{x_0}L_0$ to $T_0^*{\R}^n$ endowed with its canonical
flat connection. Let $h_{\nabla_{L_0}}$ be the holonomy
representation. Its linear part induces a representation
$h'_{\nabla_{L_0}}:\pi_1(L_0)\rightarrow T_0^*{\R}^n$. Suppose also
defined $h_{x_0}:\pi_1(L_0)\rightarrow Diff({\R}^n)$, a
representation. The homomorphism $h_{x_0}$ is the holonomy of a
lagrangian foliation which has $L_0$ as a compact leaf if and only if for every element of $D_0$ the
set of differentiable functions  $f$ defined in a neighbourhood of
$0$ in ${\R}^n$ such that $f(0)=0$ we have:
$$
h_{\nabla_{L_0}}(\gamma)\circ d_0=d_0\circ h_{x_0}(\gamma^{-1})^*,
\leqno(1)
$$
where $h_{x_0}^*(\gamma)(f)=f\circ h_{x_0}(\gamma)$.

\medskip

The representation $h'_{\nabla_{L_0}}$ endows $T^*{\R}^n$ with a
$\pi_1(L_0)$ module, we will denote by $H^*(\pi_1(L_0),T^*{\R}^n)$
the corresponding cohomology modules. The radiance obstruction map
$\gamma\rightarrow h_{\nabla_{L_0}}(0)$ defines an element
$[h'_{\nabla_{L_0}}]$ of $H^1(\pi_1(L_0),T^*{\R}^n)$.

The representation $h_{x_0}$ endows $D_0$ with a module structure
defined by
$$
\gamma\circ f=f\circ h_{x_0}(\gamma^{-1})
$$
 The map $d_0$
induces a morphism
$$
d_0^*: H^1(\pi_1(L_0),D_0)\rightarrow H^1(\pi_1(L_0),T^*_0{\R}^n)
$$

\bigskip

{\bf Theorem [2] 4.3.}

{\it The germs of lagrangian foliations in a neighbourhood of
$(L_0,\nabla_{L_0})$ which foliation holonomy is $h_{x_0}$ are
classified up to  symplectomorphisms by the elements of
${d_0^*}^{-1}[h_{\nabla_{L_0}}]$.}

\bigskip

{\bf Classification of germs of lagrangians foliation around a
compact leaf $L_0$ which Weinstein connection is not necessarily
complete.}

\bigskip

In this paragraph, we will generalize the classification of
Curras-Bosch and Molino without supposing that $(L_0,\nabla_{L_0})$
is a complete affine manifold.

\medskip

We can identify a neighbourhood of $L_0$ to a neighbourhood of the
trivial section of $T^*L_0$, we can extend the developing map
$D_0:\hat L_0\rightarrow {\R}^n$ to a map $D:T^*\hat
L_0\rightarrow {\R}^{2n}$ such that for each element $(x,y)$ of
$T^*L_0$, we have $D(x,y)=(D_0(x),y)$. This definition is possible
since $T^*\hat L_0$ is a trivial bundle.

 We consider a
 transversal $T$ diffeomorphic to a contractible open set of the
 fiber of $x_0$. The lagrangian foliation in $U_0$ is defined
 by a suspension, its lift to the universal cover $\hat U_0$ is the
 the the trivial foliation of $\hat L_0\times T$. For each element
 $x$ of $\hat L_0$, the symplectic duality allows to identify
 $T^*(x\times T)$ to ${T\hat L_0}_x$, we can also identify
 $T^*(x\times T)$  to $T^*({x_0}\times T)$. Let $U_x$ be a neighbourhood of $0$
 in $T{\hat{ L_0}}_x$ such that the restriction
of $exp_x$ associated to the affine connection to it is injective.
We can define the chart $\phi_x=exp_x: U_x\rightarrow \hat L_0$,
 the symplectic duality allows to identify $U_x$ to a an open set $V_x$ which
 contains $0$ in $T^*(x\times T)=T^*(x_0\times T)$  defines a chart
  $V_x\rightarrow L_0$.
 Following the well-known construction of the developing map,
we can define a developing map $D_0: \hat L_0\rightarrow T^*T_{x_0}.$

We have:

$$
   h_{\nabla_{L_0}}(\gamma)\circ d_{x_0}=d_{x_0}\circ h_{x_0}(\gamma^{-1})^*,
     \leqno(1)
$$
where $h_{x_0}^*(\gamma)(f)=f\circ h_{x_0}(\gamma)$.

\medskip

Conversely, let $(L_0,\nabla_{L_0})$ be a compact affine manifold,
and $h_{x_0}: \pi_1(L_0)\rightarrow Diff({\R}^n)$ a symplectic form.
Endows ${\R}^{2n}$ with a symplectic form which pulls-back by the
developing map of $T^*L_0$ is the pulls back to $T^*\hat L_0$ of
the symplectic form of $T^*L_0$. As above, we can identify using the
symplectic duality the developing
map of $(L_0,\nabla_{L_0})$ to a map $\hat L_0\rightarrow T_0^*{\R}^n$.

 The fact that
  the holonomy of $L_0$ preserve the pulls-back of this symplectic
  form to ${\R}^n\times \hat L_0$ means that:

$$
h_{\nabla_{L_0}}(\gamma)\circ d_0=d_0\circ h_{x_0}(\gamma^{-1})^*,
\leqno(1)
$$

 we thus deduce a map

 $$
d_0^*: H^1(\pi_1(L_0),D_0)\rightarrow H^1(\pi_1(L_0),T^*_0{\R}^n)
$$

We have the classification theorem

\medskip

{\bf Theorem 4.4.}

{\it The germs of lagrangian foliations in a neighbourhood of
$(L_0,\nabla_{L_0})$ which foliation holonomy is $h_{x_0}$ are
classified up to a symplectomorphism by the elements of
${d_0^*}^{-1}[h_{\nabla_{L_0}}]$.}

\bigskip

Now we classify the sequence mentioned at the beginning.
The affine manifolds $(L_i,\nabla_{L_i})$ will be supposed to be compact
and complete.
Given an isomorphism class $e_1$ of a germ of a lagrangian
foliation in a neighbourhood of $(L_1,\nabla_{L_1})$ we will
classify the germs of neighbourhoods of lagrangian foliations which
has $(L_2,\nabla_{L_2})$ has a leaf, and such that the
Marsden-Weinstein reduction of a neighbourhood $U_2$ (related to the
circle parallel action) of $L_2$ is a
symplectic manifold $(M_1,\omega_1)$ endowed with a lagrangian
foliation which has $(L_1,\nabla_{L_1})$ has compact leaf, and
which neighbourhood is isomorphic to $e_1$.

\bigskip

Let $x_2$ be an element of $L_2$ which image by the canonical
projection $p_2:L_2\rightarrow L_1$ is the element $x_1$ of $L_1$.
A germ of a lagrangian foliation on $U_2$ which has
$(L_2,\nabla_{L_2})$ has a leaf is defined by the following data:

A representation $h_{x_2}:\pi_1(L_2)\rightarrow Diff({\R}^{n+1})$,
which satisfies the relation $(1)$,

Suppose that the Marsden-Weinstein reduction of $U_2$ is a
symplectic manifold $(M_1,\omega_1)$ endowed with a lagrangian
foliation which has $L_1$ as a leaf, and which  neighbourhood is
$e_1$, then we have:

  \medskip

{\bf Proposition.}

{\it The following square is commutative:

$$
\matrix {\pi_1(L_2)&{\buildrel{h_{x_2}}\over{\longrightarrow}}&
Diff({\R}^{n+1})\cr\downarrow p_2 &\ \ \ \ & \downarrow \cr
\pi_1(L_1) &{\buildrel{h_{x_1}}\over{\longrightarrow}}&
Diff({\R}^n)} \leqno(2)
$$  }

\medskip

{\bf Proof.}

Consider the moment map $J:U_2\rightarrow {\R}$, we have supposed that
$J^{-1}(0)$ contains $L_2$. The transversal $T_2$ to the lagrangian foliation
is diffeomorphic to an open set of ${\R}^{l_1+1}$, it is enough to remark that
 we can assume that the moment map is  a coordinate function in a neighbourood
 of $0$. To do this,
 the function used to extend the action of the circle in a
 neighbourhood of $L_2$ see Molino Theorem 3.2 p. 186 must depend only of
  one coordinate in a neighbourhood
 of $0$, and be equal to this coordinate in this neighbourhood of $0$, since the
 basic function associated coincide with the moment map in this neighbourhood.

\medskip

We will denote by $e_2$ the isomorphism class of this germ. We
will endow the set $(e_1,e_2)$ where $e_1$ is fixed with a structure of a
gerbe.

\medskip

Let $Et_{L_0}$ be the site  which objects are covering spaces of
$L_0$. To each object $e$ of $Et_{L_0}$,  we define the category
$symp(e,L_0)$ of germs of lagrangian foliations which has $e$ has
a leaf, and such that the classical holonomy of an element of
$symp(e,L_0)$ in $x_e$, the element of $e$ over $x_0$ is the
restriction of $h_{x_0}$ to $\pi_1(e)$. The automorphisms of
objects of $symp(e,L_0)$ are exponentials of hamiltonians flows
defined by basic functions of the lagrangian foliation. It is easy
to check that the correspondence $e\longrightarrow symp(e,L_0)$
defines a gerbe on $Et_{L_0}$. That is the following axioms are
satisfied:

- For each map $U\rightarrow V$, between objects of $Et_{L_0}$ we
have a pulls-back map  $r_{U,V}:symp(V,L_0)\rightarrow
symp(U,L_0)$ such that $r_{U,V}\circ r_{V,W}=r_{U,W}$.

- Gluing conditions for objects,

Consider  a covering family $(U_i)_{i\in I}$  of an open set $U$ of $M$,
and for each $i$, an object $x_i$ of $symp(U_i,L_0)$. Suppose that
there exists a map $g_{ij}:r_{U_i\cap U_j,U_j}(x_j)\rightarrow
r_{U_i\cap U_j, U_i}(x_i)$ such that $g_{ij}g_{jk}=g_{ik}$, then
there exists an object $x$ of $symp(U,L_0)$ such that
$r_{U_i,U}(x)=x_i.$

Gluing conditions for arrows,

Consider two objects $P$ and $Q$ of $symp(L_0,L_0)$,  the map
$U\rightarrow Hom(r_{U,M}(P),r_{U.M}(Q))$ is a sheaf.

-There exists a covering family $(U_i)_{i\in I}$ of $Et_{L_0}$
such that for each $i$ the category $S(U_i)$ is not empty.

- Let $U$ be an object of $Et_{L_0}$, for each objects $x$ and $y$
of $symp(U,L_0)$, there exists a covering family $(U_i)_{i\in I}$
of $U$ such that $r_{U_i,U}(x)$ and $r_{U_i,U}(y)$ are isomorphic.

Every arrow of $symp(U,L_0)$ is invertible, and there exists a
sheaf $Ham$ in groups on $M$, such that for each object $x$ of
$symp(U,L_0)$, $Hom(x,x)= Ham(U)$, and the elements of this family
of isomorphisms commute with the restriction
 maps.

\medskip

{\bf Remark.}

 The gerbe $symp(L_0)$ that we have just defined is a trivial
 gerbe as shows the previous classification theorem. We will
 denote by $Symp(L_0)$.

Let $f$ be a global section of $symp(e,L_0)$, ${L_0}_f$ the pulls-back of
$L_0$ to $f$, i.e the leaf of $f$ which projects to $L_0$.
 Consider a lagrangian $1-$connected transversal $T$ of ${\cal L}_f$ to
 ${L_0}_f$,  the symplectic duality
 allows us to identify $f$ to the vertical foliation of $T_{x_0}^*T$, the
automorphisms of this objects are vertical homomorphisms of the
foliation. This group is a commutative group since $T$ is
lagrangian, and its elements are exponential of vertical hamiltonian vector
fields.

\bigskip

{\bf The general case.}

\bigskip

Consider a sequence $(L_l,\nabla_{L_l})\rightarrow...\rightarrow
(L_1,\nabla_{L_1})$ a sequence of compact and complete affine
manifolds, such that the map
$f_i:(L_{i+1},\nabla_{L_{i+1}})\rightarrow (L_i,\nabla_{L_i})$ is
an affine bundle which typical fiber is the circle, our purpose is
to classify sequences of germs, of lagrangian foliations $U_i$
such that $(L_{i+1},\nabla_{L_{i+1}})$ is a compact leaf of
$U_{i+1}$, and $L_i$ is a compact leaf of the lagrangian foliation
of the Weinstein-Marsden reduction $(M_i,\omega_i)$ of $U_{i+1}$
by the hamiltonian action of the circle.

\medskip

The germs of the lagrangian foliation of $U_i$ is defined by a
representation
$$
\pi_1(L_i)\rightarrow Diff({\R}^n)
$$
which satisfies condition $(1)$, such that the following square is
commutative:

$$
\matrix
{\pi_1(L_{i+1})&{\buildrel{h_{x_{i+1}}}\over{\longrightarrow}}&
Diff({\R}^{n+1})\cr\downarrow p_i &\ \ \ \ & \downarrow \cr
\pi_1(L_i)
&{\buildrel{h_{x_i}}\over{\longrightarrow}}&Diff({\R}^n)}
$$

Let  $e_1$ be an object of
$symp(L_1)$, we will associate to $e_1$, the gerbe
$symp(L_2,e_1)$ which is a gerbe defined on the site $Et_{L_2}$,
such that for each object $e_2$ of $Et_{L_2}$, the objects of the
category $symp(L_2,e_1)(e_2)$ are germs $U_2$ of lagrangians
foliations which has $L_2$ as a leaf, and such that
 the Marsden-Weinstein reduction of $U_2$ by the hamiltonian
 action of the circle is a manifold endowed with a lagrangian
 foliation which has $L_1$ as a  leaf and is isomorphic to $e_1$.

 \medskip

 Suppose defined the gerbe $symp(L_p,e_1,..,e_{p-1})$, and let $e_p$ be
 an object of this gerbe.
  We will
 define the gerbe $symp(L_{p+1},..,e_p)$ as the category
 which is a gerbe defined on the site $Et_{L_{p+1}}$,
 the objects of
the category $symp(L_{p+1},..,e_p)$ are germs $U_{p+1}$ of
lagrangian  foliations which has $L_{p+1}$ as
leaf, and such that
 the Marsden-Weinstein reduction of $U_{p+1}$ by the hamiltonian
 action of the circle is a manifold endowed with a lagrangian
 foliation which has $L_p$ as a  leaf, and
 such that the germ of this foliation around $L_p$ is $e_p$.

\medskip

The gerbes $symp(L_{p+1},..,e_{p})$ are trivial gerbes. For
objects $e_{p}$ and $e'_{p}$ of $symp(L_p,..,e_{p-1})$, and a
morphism $f:e_p\rightarrow e'_p$, there is a functor,
$f^*:symp(L_{p+1},..e'_p)\rightarrow symp(L_{p+1},..,e_p)$ such that
there exists an isomorphism:

$$
c(f,g):(fg)^*\rightarrow g^*f^*
$$
  which satisfies a $1-$descent condition that is:

$$
(Id*c(f,g)\circ c(fg,h) = c(g,h)*Id \circ c(f,gh)
$$

Suppose that the dimension of $L_1$ is $l_1$, then the dimension
of $L_i$ is $l_i+i-1$. We will denote by $D_i$ the germs at $0$ of
differentiable functions of ${\R}^{l_i+i-1}$ which take the value
$0$ at the origin. The holonomy in $x_i$ endows $D_i$ with the
structure of a $\pi_1(L_i)$ module by setting:
$$
\gamma\circ f=f\circ h_{x_i}(\gamma^{-1})
$$
the differential at the origin $d_i$ induces a morphism:
$$
d_i^*:H^1(\pi_1(L_i),D_i)\longrightarrow
H^1(\pi_1(L_i),T^*_0{\R}^{l_1+i-1})
$$
The isomorphisms classes of germs of lagrangian foliation for
which the holonomy is $h_{x_i}$ are the elements of
${d_i^*}^{-1}[h_{\nabla_{L_i}}]$.

Given a global section $e_i$ of $symp(L_i)$, classified by
the class $[c_i]$ of $H^1(\pi_1(L_i),D_i)$, we will determine the
the classifying cocycle of the global sections of $symp(L_{i+1})$ which
gives rise to $e_i$.

The surjection $p_i:\pi_1(L_{i+1})\rightarrow \pi_1(L_i)$, and the
commutative square $(1)$ induces a the following commutative
square:

$$
\matrix{H^1(\pi_1(L_{i+1}),D_{i+1})&\longrightarrow&
H^1(\pi_1(L_i),{\R}^{l_1+i})\cr f_i\downarrow && g_i\downarrow \cr
H^1(\pi_1(L_i),D_i)&\longrightarrow&
H^1(\pi_1(L_i),{\R}^{l_1+i-1})}.
$$
The image of the classifying cocycle of elements of
$symp(L_{i+1})$ by $f_i$ is the classifying cocycle of
$e_i$.

\medskip

Let $\gamma_1,...,\gamma_{l_1+i-1}$ be the generators of
$\pi_1(L_i)$, we will suppose that if $j\leq i$,
$\gamma_1,...,\gamma_{l_1+j-1}$ are generators of $\pi_1(L_j)$.
Remark that since we have supposed that $(L_2,\nabla_{L_2})$ is a
circle bundle over $(L_1,\nabla_{L_1})$, we can assume that
$h_{x_2}(\gamma_{l_1+1})$ is the identity, that is the generator
of the element of $\pi_1(L_2)$ which preserves leaves of the
circle bundle.

Let $c_1$ be the cocycle which defines the isomorphism class of
the lagrangian foliation ${\cal L}_1$ on $U_1$,  for each element
$\gamma$ of $\pi_1(L_1)$, $c_1(\gamma)$ is a germ of a
differentiable function defined on ${\R}^{l_1}$. Let $\gamma'$ be
an element of $\pi_1(L_2)$ over $\gamma$, then $c_2(\gamma')$ is a
germ of a differentiable function of ${\R}^{l_1+1}$ over
$c_1(\gamma)$, where $c_2$ is a classifying cocycle of an element
of $symp(L_2,e_1)$.

\medskip

Consider $L(h_i)$, the linear holonomy of the affine manifold
$(L_i,\nabla_{L_i})$, and write ${\R}^{l_1+1}={\R}^{l_1}+{\R}$.
For each element $\gamma\in\pi_1(L_2)$, $L(h_2(\gamma))$ depends of
 $p_2(\gamma)$ and the
 projection of $L(h_2)(\gamma)$ on
${\R}e_{l_1+1}$ parallel to ${\R}^{l_1}$ is a cocycle $d_2$ for
the trivial action of $\pi_1(L_2)$ on ${\R}e_{l_1+1}$.

\medskip

We will suppose now that the holonomy is linearizable in  neighbourhoods
compact leaves $L_1$ and $L_2$. We will determine the relation between the
cocycle which define elements of $symp(L_2,e_1)$ and the classifying cocycle
of $e_1$.

For each element $\gamma$ of $\pi_1(L_1)$, we consider a differentiable
function $c_2(\gamma)$ of ${\R}^{l_1+1}$ such that $c_2(\gamma)$ project to
$c_1(\gamma)$, which
means that there exists a function $f_2(\gamma)$ such that

$$
c_2(\gamma)=c_1(\gamma)(x_1,..,x_{l_1})+f_2(\gamma)(x_{l_{1+1}})
$$

 and for $\gamma_{l_1+1}$, we consider an element
$c_2(\gamma_{l_1+1})$ which projects
to $0$, that is a function of $x_{l_1+1}$ the $l_1+1$ coordinate in the basis
$(e_1,...,e_{l_1+1})$ of ${\R}^{l_1+1}$, the fact that the map $c_2$ is a
cocycle implies the following:

$$
c_2(\gamma_{l_1+1}\gamma)=c_2(\gamma_{l_1+1})+c_2(\gamma)
$$
For every element $\gamma$ of $\pi_1(L_2)$,

$$
c_2(\gamma\gamma')=c_2(\gamma)+\gamma c_2(\gamma')=
$$

$$
c_1(p_2(\gamma))+f_2(\gamma)+ p_2(\gamma)c_1(p_2(\gamma)) +
f_2(\gamma')\circ d_2(\gamma)
$$

This implies that:

$$
f_2(\gamma\gamma')=f_2(\gamma)+f_2(\gamma')\circ d_2(\gamma)
$$

More generally,  the classifying cocycle $c_i(h)$ of an element
 $h$  of $symp(L_i)$
is described by a germ of a differentiable function $f_i(h)$ of
${\R}e_{l_1+i-1}$, such that

$$
c_i(h)(\gamma)= c_{i-1}(p_i(\gamma)+f_i(h)(\gamma)
$$

 which satisfies:

$$
 f_i(h)(\gamma\gamma')= f_i(h)(\gamma)+
 f_i(h)(\gamma)\circ d_i(\gamma).
$$

Cocycles in this theory are computable when we suppose that the
lagrangian foliation is linearizable.

\bigskip

{\bf Bibliography.}

\medskip

[1] Charette, V. Goldman, W.

  Affine schottky groups and Crooked tilings. Proceedings of Crystallographic
  groups and their generalizations II. Kortrijk 1999.

[2] Curras-Bosch, C. Molino, P.

Holonomie suspension et
classification pour les feuilletages lagrangiens. C.R.A.S (326)
1317-1320

[3] Dazord, P.

Sur la geometrie des sous-fibres et des feuilletages lagrangiens.
Annales scient. E. Norm. Sup. Paris.

[4] Fried, D. Goldman, W.

Three-dimensional affine
 crystallographic groups. Advances in Math. 47 (1983), 1-49.

[5] Fried, D. Goldman, W. Hirsch, M.

 Affine manifolds with
 nilpotent holonomy. Comment. Math. Helv. 56 (1981) 487-523.

 [6] Godbillon, C.

 Feuilletages. Etudes g\'eom\'etriques. Progress
 in Mathematics, 98.

 [7] Margulis, G.

 Complete affine locally flat manifolds with a
 free fundamental group. J. Soviet. Math. 134 (1987), 129-134.

 [8] Milnor, J. W.

 On fundamental groups of complete affinely flat
 manifolds, Advances in Math. 25 (1977) 178-187.

[9] Molino, P Lagrangian holonomy.

Analysis and geometry in foliated manifolds. (Santiago de Compostela) 183-194.

[10] Plante, J.F. Polycyclic groups and transversely affine foliations.

J. Diff. Geometry 35 (1992) 521-534.

[11] Plante, J. F. Foliations with measures preserves holonomy.

Ann. of Math. 102 (1975) 327-361.

[12] Plante J.F Thurston, W. Polynomial growth in holonomy foliations.

Comment. Math. Helv. 51 (1976) 567-584.

[13] Sullivan, D. Thurston, W. Manifolds with canonical

 coordinate charts: some examples. Enseign. Math 29 (1983), 15-25.

 [14] Tsemo, A. Dynamique des vari\'et\'es affines. J. London Math.
 Soc. 63 (2001) 469-487.

[15] Tsemo, A. Non abelian cohomology: The viewpoint of gerbed towers.
Preprint.

[16] Weinstein, A. Symplectic manifolds and their lagrangian submanifolds.

Advances in Math 6 (1971) 329-346.

\end{document}